\documentclass{amsart}
\usepackage{url}
\usepackage{nohyperref}
\begin{document}
\noindent
Topology Atlas Invited Contributions \textbf{8} (2003) no.~3, 1--2.
\smallskip
\title{Supernearness spaces as a tool for studying Unification and extensions}
\author{D. Leseberg}
\address{Dieter Leseberg\\
Department of Mathematics and Informatics\\
Free University of Berlin\\
Germany}
\email{d.leseberg@tu-bs.de}
\date{October 2003}
\maketitle

Topological extensions are closely related to nearness structures of
various kinds. For example, the Smirnov compactification of a proximity
space $X$ is a compact Hausdorff space $Y$ which contains $X$ as a dense
subspace and for which it is true that a pair of subsets of $X$ is near
if, and only if, their closures in $Y$ meet.  Lodato generalized this
result to weaker conditions for the proximity and the space using
\emph{bunches} for the characterization of the extension.

Ivanova and Ivanov studied contiguity spaces and bicompact extensions such
that a finite family of subsets are contigual if, and only if, there is a
point of $Y$ which is simultaneously in the closure in $Y$ of each set of
the family. Herrlich found a useful generalization of contiguity spaces by
introducing nearness spaces, and Bentley showed that those nearness spaces
which can be extended to topological ones have a neat internal
characterization. Doitchinov introduced the notion of supertopological
spaces in order to construct a unified theory of topological, proximity
and uniform spaces, and he proved a certain relationship of some special
classes of supertopologies---called $b$-supertopologies---with compactly
determined extensions. Recently, supernear spaces were introduced by
myself in order to define a common generalization of nearness spaces and
supertopological spaces as well. As a basic concept, the notion of bounded
sets is used, and functions (supernear operators) from a collection of
bounded sets into powers of a given set are considered which fulfill
certain axiom of nearness. So, we define for each bounded set the
so-called \emph{$B$-near collections}, and a bounded map is defined as a
supernear-map if it preserves such collections. Thus, the new topological
category $\mathbf{SN}$ whose objects are the supernear spaces is
established, and moreover we can characterize those supernear spaces which
can be extended to topological ones. As a basic notion we consider
topological extensions which are strict in the sense of Banaschewski,
where the hull in $Y$ of the images of sets in $X$ form a base for the
closed subsets in $Y$. In this connection we study symmetrical and
non-symmetrical extensions, and consequently we can characterize and
describe those supernear spaces which are in a certain one-to-one
relationship or correspondence, respectively, with them.

In other words, we describe those supernear spaces, whose $B$-near
collections are determined by their respective constructing strict
topological extension. As special cases we get back the former described
constructions by various authors. So it seems to be possible to describe
topological unification and extensions by the support of one 
\emph{topological concept} with its corresponding maps.

\providecommand{\bysame}{\leavevmode\hbox to3em{\hrulefill}\thinspace}

\end{document}